\input amstex
\documentstyle{amsppt}
\magnification=1200
\hcorrection{.25in}
\advance\voffset.25in
\advance\vsize-.5in

\topmatter
\title Canonical systems and
finite rank perturbations of spectra
\endtitle
\author Alexei G.  Poltoratski\endauthor
\thanks Research at MSRI is supported in part
by NSF grant DMS-9022140.\endthanks
\address MSRI, 1000 Centennial Drive, Berkeley, California 94720\endaddress
\email agp\@msri.org \endemail
\abstract We use Rokhlin's Theorem on the uniqueness of canonical systems
to find a new way to establish connections between Function Theory
in the unit disk and rank one perturbations of self-adjoint or unitary
operators. In the n-dimensional case, we prove that
for any cyclic self-adjoint operator $A$,
operator
$A_\lambda= A + \Sigma_{k=1}^n \lambda_k(\cdot,\phi_k)\phi_k$
is pure point for a.e.\
$\lambda=(\lambda_1,\lambda_2,...,\lambda_n)
\in\Bbb R^n$ iff 
operator
$A_\eta=A+\eta(\cdot,\phi_k)\phi_k$
is pure point for a.e.\ $\eta\in\Bbb R$ for $k=1,2,...,n$.
We also show that if $A_\lambda$ is pure point for a.e.\ $\lambda\in
\Bbb R^n$ then $A_\lambda$ is pure point for a.e.\ $\lambda\in \gamma$
for any analytic curve $\gamma\in\Bbb R^n$.
\endabstract
\tolerance=1000

\endtopmatter
\document

\heading Introduction \endheading

This note analyzes the spectral properties of finite rank perturbations
of self-adjoint and unitary operators. In Section 1
we will try to put some well-known results on rank one perturbations
into a
different prospective.
The new approach will allow us to shorten some of the existing proofs.
Section 2 deals with rank n perturbations and contains the main results
of this paper.

Our approach in Section 1
will be based on the notion of canonical systems of measures.
This notion  was introduced by
Kolmogorov and Rokhlin [R], and, independently, by Ambrose, Halmos
and Kakutani [A-H-K, H1-2] more
than half a century ago. It was originally used to study partitions,
homeomorphisms and factor spaces of various measure spaces.
The term ``canonical'' belongs, most likely, to Kolmogorov.
We will first discuss examples of canonical systems
of measures naturally appearing in certain problems of Function
Theory in the unit disk and Perturbation Theory. We will then
use Rokhlin's uniqueness Theorem to show that different
problems yield the same families of measures.

In our first example of canonical systems (see Example 1 in Section 1)
we are going to discuss families of measures on the unit circle
generated by inner functions in the unit disk. Such families
first appeared in [C] and then were further studied
in [A1-5], [P1-3] and [S]. Our second example (see Example 2 in Section 1)
deals with families of spectral measures of rank one perturbations
of self-adjoint or unitary operators. These families were
extensively studied by many authors. We refer to [A], [A-D] and [D] for
the basic results on this subject, and to  [J-L-R-S], [R-M-S], [P3] and [Si]
for the latest developments
and further references.

As was proven by D. Clark in [C], any family of measures generated by
an inner function (as in Example 1) is a family of spectral measures
of all rank one perturbations of a model contraction. In Section 1
we will show how Rokhlin's Theorem can help establish this connection
for families of spectral measures of cyclic self-adjoint or unitary
operators.

Section 2 is devoted to n-dimensional analogies of Simon-Wolff
Theorem on the pure point spectrum of a random rank one perturbation.
The one-dimensional result was applied in [S-W] to prove the
the existence of Anderson localization for the  wave propagation
described by the discrete Schr\"odinger operator with random potential
in one dimension.
It is still an open question if Anderson localization takes place
in dimension more than one.
 
As we will find out, the n-dimensional case adds certain new and
somewhat unexpected features to the general picture. 
First, using the observations made in Section 1, 
we will prove that
for any cyclic self-adjoint operator $A$ and its cyclic vectors
$\phi_1,...,\phi_n$,
operator
$A_\lambda=A+\sum_{k=1}^n \lambda_k(\cdot,\phi_k)\phi_k$
is pure point for a.e.\
$\lambda=(\lambda_1,\lambda_2,...,\lambda_n)
\in\Bbb R^n$ iff
operator
$A_\eta=A+\eta(\cdot,\phi_k)\phi_k$
is pure point for a.e.\ $\eta\in\Bbb R$ for $k=1,2,...,n$.
I. e.  random perturbation
$A_\lambda$ is pure point almost surely iff
it is pure point almost surely on the skeleton $S$ (the
union of coordinate axes) in $\Bbb R^n$.
Using the terminology of Function Theory,
this result means that $S$ is, in some sense, a {\it sampling set} for
this perturbation problem.

We will also prove that
if $A_\lambda$ is pure point for a.e.\ $\lambda\in
\Bbb R^n$ then $A_\lambda$ is pure point for a.e.\ $\lambda\in \gamma$
for any analytic curve $\gamma\in\Bbb R^n$.

This statement reveals a certain fine structure of the set
$$E=\{\lambda\in\Bbb R^n|
A_\lambda \text {is not pure point }\}.$$
Such a property of $E$ may seem surprising, 
since various examples of perturbations of singular spectra
show that $E$ can be almost  ``arbitrarily bad''.
For instance, as follows from the results of [R-M-S] and [G],
if the spectrum of $A$ contains an interval, then,
even if $A_\lambda$ is pure point a. s., $E$
is everywhere dense and $G_\delta$
on every line parallel to the coordinate axes.

The results of Section 2 allow some infinite-dimensional generalizations
which are to be published elsewhere.

\remark{\bf Acknowledgments} This work was done during my stay at
the Mathematical
Sciences Research Institute, and
I would like to thank its administration
and staff for their hospitality and support.
\endremark

\heading 1. Canonical systems of measures in function theory
and functional analysis \endheading

Let us start with the following

\definition{Definition}
Consider a measure space $(X, \Cal A, \mu)$ where $X$ is a set of points,
$\Cal A$ is a $\sigma$-algebra of subsets of $X$ and $\mu$ is a finite measure
on $\Cal A$. Space $(X, \Cal A,\mu)$ is called Lebesgue space if there exist a
1-1 measure preserving map from $(X, \Cal A, \mu_c)$ onto the unit
interval with Lebesgue measure.
\enddefinition

Here $\mu_c$ denotes the continuous part of $\mu$ and the term
"measure preserving" means in particular that the map induces
a 1-1 correspondence between $\Cal A$ and Lebesgue $\sigma$-algebra on the unit
interval.

Let $\xi$ be a measurable partition of a Lebesgue space $(X, \Cal A, \mu)$:
$\xi=\{S_\alpha\}_{\alpha\in Y}$ where $Y$ is some set of parameters,
$S_\alpha\in\Cal A$, $S_\alpha\cap S_\beta=\emptyset$ if $\alpha\neq\beta$
and $\bigcup_Y S_\alpha=X$. Then we can consider the factor space
$(X/\xi, \Cal A_\xi, \mu_\xi)$ where $\sigma$-algebra $\Cal A_\xi$ and measure 
$\mu_\xi$ are induced
by $\Cal A$ and $\mu$ via the factor map.

Suppose each set $S_\alpha$ from our partition $\xi$ is itself
a Lebesgue space $(S_\alpha,\Cal A_\alpha, \mu_\alpha)$ for some $\sigma$-algebra
$\Cal A_\alpha$ and measure $\mu_\alpha$.

\definition{Definition} The system of measures $\{\mu_\alpha\}_{\alpha\in Y}$
 is called
canonical if for any $A\in\Cal A$ we have
$A\cap S_\alpha\in\Cal A_\alpha$ for $\mu_\xi$-a.e.\
$\alpha$ and 
$$\mu(A)=\int_{X/\xi}\mu_\alpha(A\cap S_\alpha) d\mu_\xi(\alpha).   \tag 1$$
\enddefinition

Note that the factor space $X/\xi$ can be naturally identified with the
set $Y$ of parameters $\alpha$. This justifies the
expression ``$\mu_\xi$-a.e.\ $\alpha$'' and integration over
$d\mu_\xi(\alpha)$ in (1).

Suppose we have a measurable partition $\xi$ of $(X,\Cal A,\mu)$.
Does there exist
a corresponding canonical system of measures? If $\xi$ is countable
the answer is obviously yes. In this case
 $\mu_\xi$ in (1) is discrete and $\mu_\alpha$ can be chosen as
the restriction of $\mu$ on $S_\alpha$. 
However when $\xi$ is uncountable, most of $S_\alpha$ have measure
0 and the canonical system can not be constructed so easily.
Nevertheless we have the following

\proclaim{Theorem 1 (Rokhlin [R])}
For any measurable partition $\xi$ of Lebesgue measure space
$(X, \Cal A, \mu)$ there exists a canonical
system $\{\mu_\alpha\}$ satisfying (1). Such a system
is essentially unique i. e. for any  two such canonical systems
$\{\mu_\alpha\}$ and $\{\mu'_\alpha\}\ \ \ $  $\mu_\alpha=\mu'_\alpha$ for $\mu_\xi$-a.e.\ $\alpha$.
\endproclaim

In the rest of this section we are going to use the uniqueness part
of this result to show how different problems in function theory
and functional analysis yield the same families of measures.

Note that condition $\bigcup S_\alpha=X$
can be replaced by a weaker condition
$$\mu(\bigcup S_\alpha)=\mu(X)\tag 2$$
in all the above definitions and in Theorem 1.

Our first example of a canonical system is related to  Analytic
Function Theory in the unit disk.

\example{Example 1}
Let $\theta$ be an inner function in the unit disk $\Bbb D$.
We will denote by $M_+(\Bbb T)$ the set of all finite positive Borel
measures on the unit circle $\Bbb T=\partial \Bbb D$.

For each $\alpha\in\Bbb T$ function $\frac{\alpha+\theta}{\alpha-\theta}$
has positive real part in $\Bbb D$. Therefore for each $\alpha\in\Bbb T$
there exists $\mu_\alpha\in M_+(\Bbb T)$ such that its Poisson integral
satisfies
$$P\mu_\alpha=\int_{\Bbb T}\frac{1-|z|^2}{|\xi-z|^2} d\mu_\alpha(\xi)
=Re\frac{\alpha+\theta}{\alpha-\theta}
.\tag 3$$
We will denote by $M_\theta(\Bbb T)$ the family $\{\mu_\alpha\}_{\alpha\in\Bbb T}$
of all such measures corresponding to $\theta$. Note that since
$Re\frac{\alpha+\theta}{\alpha-\theta}=0$ a.e.\ on $\Bbb T$,
all measures $\mu_\alpha$ are singular.

If $\theta$ is defined in the upper half plane $\Bbb C_+$ then one can replace
the Poisson integral in (3) with the Poisson integral in $\Bbb C_+$ and define
an analogous family $M_\theta(\hat\Bbb R)$ on the completed real line
$\hat\Bbb R=\Bbb R\cup\{\infty\}$ (see [P3]).

Families $M_\theta$ possess many interesting properties (see
[A1-5], [C], [P1-3] and [Sa]). Our next task is to show that
all such families are {\it canonical systems}.

First let us notice that the defining formula (3) implies that
measures $\mu_\alpha$ are concentrated on disjoint sets
$S_\alpha=\{\xi\in\Bbb T | \theta(\xi)=\alpha\}$ (all analytic functions
in this paper are defined on the boundary by their nontangential boundary
values). Also if $\Cal P_z$ is a Poisson kernel for some $z\in\Bbb D$ then
$$\int_{\Bbb T}\int_{\Bbb T} \Cal P_zd\mu_\alpha dm(\alpha)=\int_{\Bbb T}
Re\frac{\alpha+\theta(z)}{\alpha-\theta(z)}dm(\alpha)=$$
$$=1=\int_{\Bbb T} \Cal P_z dm.\tag 4$$
Note that since linear combinations of Poisson kernels are dense
in the space of all continuous functions on $\Bbb T$, one can replace
$\Cal P_z$ in (4) with an arbitrary continuous function. Further,
by taking the limit over continuous functions, one can replace $\Cal P_z$
 in (4)
with a characteristic function of any Borel set $B\subset \Bbb T$ to
obtain the following analogy of formula (1) (see [A1]):
$$\int_{\Bbb T} \mu_\alpha(B)dm(\alpha)=m(B).\tag 5$$

Therefore family $M_\theta$ defined by formula (3) is indeed a
canonical system on $(\Bbb T, \Cal B, m)$ for the partition
$S_\alpha=\{\theta=\alpha\}$ where $\Cal B$ is Borel
$\sigma$-algebra on $\Bbb T$.

The same argument holds in the case of the upper half plane for the families
$M_\theta(\hat\Bbb R)$.
\endexample

Our second example concerns Perturbation Theory of self-adjoint
and unitary operators.

\example{Example 2} 
Let $A_0$ be a singular bounded cyclic self-adjoint operator
acting on a separable Hilbert space, $\phi$
its cyclic vector. Consider the family of rank one perturbations of $A_0$:
$$A_\lambda=A_0+\lambda(\cdot,\phi)\phi\ \ \ \ \ \lambda\in\Bbb R.$$
Let $\mu_\lambda$ be the spectral measure of $\phi$ for $A_\lambda$.
We will denote by $M_{A_0,\phi}$ 
the family $\{\mu_\lambda\}_{\lambda\in\Bbb R}$.

For basic results and references on families $M_{A_0,\phi}$ the reader can
consult [A], [A-D], [D] or [S].

Here we will show that family $M_{A_0,\phi}=\{\mu_\lambda\}_{\lambda\in\Bbb R}$
is a {\it canonical system}.

Note that Cauchy integral of $\mu_\lambda$ satisfies:
$$K\mu_\lambda(z)=\int_{-\infty}^{\infty}\frac{d\mu_\lambda(t)}
{t-z}=((A_\lambda-z)^{-1}\phi,\phi).\tag 6$$
Also for the resolvents of operators $A_\lambda$ we have:
$$(A_0-z)^{-1}-(A_\lambda-z)^{-1}=\left[\lambda(\cdot,\phi)
(A_\lambda-z)^{-1}\phi\right](A_0-z)^{-1}.\tag 7$$
Combining (6) and (7) we obtain
$$K\mu_\lambda(z)=\frac{K\mu_0(z)}{1+\lambda K\mu_0(z)}\tag 8$$
(see [A]).

Let us observe that, as follows from (6), measures $\mu_\lambda$ are
concentrated on the sets $S_\lambda=\{x\in\Bbb R| ((A_\lambda-x)^{-1}\phi,\phi)=
\infty\}$.
Via (7) $S_\lambda$ can be redefined as
$$S_\lambda=
\{x\in\Bbb R| ((A_0-x)^{-1}\phi,\phi)=
-\frac1\lambda\}$$
 which shows that $S_\lambda$ are pairwise disjoint.

To prove the integral formula notice that
if $\Cal P_z=\frac y{(x-t)^2+y^2}$ is the Poisson kernel for $z=x+iy\in\Bbb C_+$
then
$$\int_{-\infty}^\infty\int_{-\infty}^{\infty}\Cal P_z(t)d\mu_\lambda(t)d\lambda=
\int_{-\infty}^\infty Re \frac{K\mu_0(z)}{1+\lambda K\mu_0(z)}d\lambda.$$
Since $\mu_\lambda$ is a positive measure, $\omega=-1/(K\mu_0(z))$
belongs to $\Bbb C_+$. Therefore
$$\int_{-\infty}^\infty\int_{-\infty}^{\infty}\Cal P_z(t)d\mu_\lambda(t)d\lambda=
\int_{-\infty}^{\infty}\Cal P_\omega(\lambda)d\lambda
=\pi=
\int_{-\infty}^{\infty}\Cal P_z(\lambda)d\lambda.$$
Since every continuous $L^1(\Bbb R)$ function can be approximated by linear
combinations of Poisson kernels, in the same way as in the previous example
we can conclude that
$$\int_{-\infty}^{\infty}\mu_\lambda(B)d\lambda=|B|\tag 9$$
for any Borel set $B$ (see [S]).

Therefore family $M_{A_0,\phi}$ is a canonical system on $(\Bbb R, \Cal B, dx)$
for the partition $S_\lambda=\{((A_0-x)^{-1}\phi,\phi)=-\frac 1\lambda\}$,
where $\Cal B$ denotes Borel $\sigma$-algebra.

An analogous argument can be applied in the case of unitary perturbations.

Let $U_1$ be a unitary cyclic singular operator, $v, ||v||=1$
its cyclic vector.
Consider the family of
unitary rank one perturbations of $U_1$ : $$U_\alpha=
U_1+(\alpha-1)(\cdot , U_1^{-1}v)v, $$ $\alpha\in\Bbb T$.
For the resolvents we have
$$(U_1-z)^{-1}-(U_\alpha-z)^{-1}=
(U_\alpha-z)^{-1}[(\alpha-1)(\cdot , U_1^{-1}v)v]
(U_1-z)^{-1}
.$$ 

Denote by $\nu_\alpha$ the spectral measure of $v$ for
$U_\alpha$. Let $ K\nu_\alpha$  be the  Cauchy 
integral of measure $\nu_\alpha$ in the unit disk $\Bbb D$:
$$K\nu_\alpha=\int_{\Bbb T}\frac{1}{1-\overline\xi z}d\nu_\alpha(\xi).$$
Since
$$ K\nu_\alpha=((U_\alpha-z)^{-1}v , U_\alpha^{-1}v)$$
for $z\in\Bbb D$,
we have that
$$ K\nu_\alpha=\frac{\alpha K\nu_1}
              {1+ (\alpha-1) K\nu_1} .\tag 10 $$
Using (10) in the same way as we  used (7) we can obtain that family
$M_{U_1,v}=\{\nu_\alpha\}_{\alpha\in\Bbb T}$ is a canonical system on
$(\Bbb T, \Cal B, m)$ for the partition 
$$ S_\alpha=\{((U_1-z)^{-1}v , U_1^{-1}v)=\frac 1{1-\alpha}\}.$$
\endexample

Now we will use Rokhlin's Theorem to show that
every family of spectral measures of rank one
perturbations (Example 2) is in fact a family
generated by an inner function (Example 1)  and vice versa.

Indeed, if we consider the following inner function in $\Bbb C_+$:
$$\theta(z)=\frac {1-i((A_0-z)^{-1}\phi,\phi)}{1+i((A_0-z)^{-1}\phi,\phi)}
\tag 11$$
then the partitions from Examples 1 and 2 coincide: $S_\alpha=S_\lambda$
for $\alpha=\frac {\lambda+i}{\lambda-i}$. Hence by Rokhlin Theorem
canonical system $M_\theta$ essentially coincides with canonical
system $M_{A_0,\phi}$. But in both systems the measures depend on
parameter continuously (with respect to the $*$-weak
topology of the space of measures). Since
the systems coincide {\it almost everywhere} and are continuous,
they must coincide {\it everywhere}:
$$M_\theta=M_{A_0,\phi}.\tag 12$$

Conversely, for each inner $\theta$ we can  choose $A_0$ to be the operator
of multiplication by $z$  in $L^2(\mu_1), \ \ \mu_1\in M_\theta$ and
$\phi=1\in L^2(\mu_1)$. Then $\theta$, $A_0$, and $\phi$ will satisfy (12).

In a similar  way, for unitary operators, if we choose
$$\theta(z)=\frac 1{K\nu_1(z)}-1$$
for any $z\in\Bbb D$ then we will have
$$M_\theta=M_{U_1,v}.\tag 13$$
with $\mu_\alpha=\nu_\alpha$ for any $\alpha\in\Bbb T$, $\mu_\alpha\in M_\theta
(\Bbb T)$ and $\nu_\alpha\in M_{U_1,v}$
(cf. [C]).

The families of measures from Examples 1 and 2 are not the only natural
examples of canonical systems. For instance, in [D] one can find two
more families of operators: the family of all self-adjoint extensions
of a given symmetric operator with deficiency indices (1,1) and
the family of all self-adjoint extensions associated with a given
limit-point Stourm-Liouville problem. Applying the argument, similar to
the one we used here, one can prove that the corresponding
families of spectral measures are canonical systems generated by inner
functions.

If $A_0$ and $\phi$ are as in Example 2, then $A_i=A_0+i(\cdot,\phi)\phi$
is a dissipative operator. Its Cayley transform $T$ is a completely
nonunitary $C_0$ contraction. Sz.Nagy-Foias model theory
provides $T$ with a characteristic
function $I$ in the unit disk. Function $I$  is related to our function $\theta$,
defined by (11), via the equation $I(z)=\theta(\omega(z))$ where $\omega$
is the standard conformal map from the unit disk onto the upper half-plane.

\heading 2. Finite rank perturbations of singular spectra \endheading

In this section we will need the
following Lemma, which can be obtained directly from the definition
of families $M_\theta$:

\proclaim{Lemma 2}
Let $\theta$ be an inner function in $\Bbb D$ and $\xi\in\Bbb T$. Then

1) measure $\mu_\alpha\in M_\theta$ has a point mass at $\xi$ iff
functions $\theta$ and $\theta'$ have non-tangential limits $\theta(\xi)$
and $\theta'(\xi)$ at $\xi$
and $\theta(\xi)=\alpha$;

2) function $\theta'$ has a non-tangential limit $\theta'(\xi)$ at $\xi$
iff there exists $\mu_\beta\in M_\theta$ such that
$$\int_\Bbb T\frac{d\mu_\beta(\omega)}{|\xi-\omega|^2}<\infty;\tag 14$$

3) measure $\mu_\alpha\in M_\theta$ has a point mass at $\xi$ iff
the non-tangential limit $\theta(\xi)$ is equal to $\alpha$ and
(14) holds for any $\beta\in\Bbb T,\ \beta\neq\alpha$.

\endproclaim

An analogous statement holds true for families $M_\theta$ on the real line.

Suppose $\theta'(\xi)$ (the  non-tangential limit of $\theta'$ at $\xi$)
exists
everywhere on $\Bbb T$ except for $\xi\in E, \ E\subset\Bbb R$.
If $m(E)=0$ then,
by formula (5), for measures from $M_\theta$ we have
$\mu_\alpha(E)=0$ for almost every $\alpha\in\Bbb T$. Since $\mu_\alpha$
is concentrated on $\{\theta=\alpha\}$, statement 1) from Lemma 2 implies
that $\mu_\alpha$ is pure point for a.e.\ $\alpha$. Now we can use
statement 2) from Lemma 2 and connections with operator theory discussed
in the previous section (formula (12)) to obtain the following

\proclaim{Theorem 3 [S-W]}
Let $A$ be a cyclic self-adjoint operator, $\phi$ its cyclic vector and
$\mu$ the spectral measure of $\phi$ for $A$. Then operator
$$A_\lambda=A+\lambda(\cdot,\phi)\phi$$
is pure point for a.e.\ $\lambda\in\Bbb R$ iff
$$\int_\Bbb R\frac{d\mu(x)}{|x-y|^2}<\infty\tag 15$$
for a.e.\ $y\in\Bbb R$.
\endproclaim

The purpose of this section is to find an n-dimensional analogy
of this result.

We will prove the following

\proclaim{Theorem 4}
Let $A$ be a cyclic self-adjoint operator, $\phi_1,\phi_2,...,\phi_n$
its cyclic vectors.

The following two conditions are equivalent:

1) operator
$$A_\lambda=A+\sum_{k=1}^n \lambda_k(\cdot,\phi_k)\phi_k$$
is pure point for a.e.\
$\lambda=(\lambda_1,\lambda_2,...,\lambda_n)
\in\Bbb R^n$.

2) operator
$$A_\lambda=A+\lambda(\cdot,\phi_k)\phi_k$$
is pure point for a.e.\ $\lambda\in\Bbb R$ for $k=1,2,...,n$;

i. e.
$$\int_\Bbb R\frac{d\mu_k(x)}{|x-y|^2}<\infty$$
for a.e.\ $y\in\Bbb R$ for $k=1,2,...,n$,
where $\mu_k$ is the spectral measure of $\phi_k$ for $A$.

\endproclaim

The theory of unitary rank n perturbations is slightly more complicated
due to the fact that a sum of two unitary perturbations may not be unitary.
Because of that, we will have to define the rank n perturbation
of a  unitary operator $U$ corresponding to vectors $\phi_1,...,\phi_n$
recursively.
I. e., let $\alpha=(\alpha_1,...,\alpha_n)\in\Bbb T^n$. Denote
$$\alpha^1=\alpha_1\in\Bbb R,\  \alpha^2=(\alpha_1,\alpha_2)\in\Bbb R^2,...$$
$$...,\alpha^{n-1}=(\alpha_1,...,\alpha_{n-1})\in\Bbb R^{n-1},
\ \ \text {and}\ \ \alpha^{n}=\alpha=(\alpha_1,...,\alpha_{n})\in\Bbb R^{n}.$$
Put
$$U_{\alpha^1}=U+(\alpha_1-1)(\cdot,U^{-1}\phi_1)\phi_1,$$
and
$$U_{\alpha^k}=U_{\alpha^{k-1}}+(\alpha_k-1)(\cdot,U_{\alpha^{k-1}}
^{-1}\phi_k)\phi_k$$
for $k=2,3,...,n$. Finally, denote $U_\alpha=U_{\alpha^n}$.

If vectors $\phi_1,\phi_2,...,\phi_n$ are pairwise orthogonal,
then our definition yields
$$U_\alpha=U+\sum_{k=1}^n\alpha_k(\cdot,U^{-1}\phi_k)\phi_k.$$

One can show that the following Theorem is
equivalent to Theorem 4:

\proclaim{Theorem 4'}
Let $U$ be a cyclic unitary operator, $\phi_1,\phi_2,...,\phi_n$
its cyclic vectors.

Then the following two conditions are equivalent:

1) operator
$U_\alpha$ (defined as above)
is pure point for a.e.\ $\alpha\in\Bbb T^n$;

2) operator
$U_\alpha=U+\alpha(\cdot,U^{-1}\phi_k)\phi_k$
is pure point for a.e.\ $\alpha\in\Bbb T$
for $k=1,2,...,n$.

\endproclaim

Since it is more convenient for us to operate in the unit
disk,
we will prove the result in the form of Theorem 4'
rather than Theorem 4. Our proof will be based on the following approach.

Let $U_1$ be a cyclic singular unitary operator. Suppose $\phi_1$ and
$\phi_2$ are its cyclic vectors, $||\phi_1||=||\phi_2||=1$.
Let $\mu_\alpha$ be the spectral measure
of $\phi_1$ for
$$U_\alpha=U_1+(\alpha-1)(\cdot,U_1^{-1}\phi_1)\phi_1.\tag 16$$
Then, as was shown in Section 1,
there exists an inner function $\theta,
\ \theta(0)=0$ such that $\{\mu_\alpha\}_{\alpha\in\Bbb T}=M_\theta$.
By the Spectral Theorem, each operator $U_\alpha$ is unitarily
equivalent to the operator $Y_\alpha$ of multiplication by $z$ in
$L^2(\mu_\alpha)$.

We will denote by $\theta^*(H^p), p>0$
the invariant subspace of the backward
shift operator in $H^p$ corresponding to $\theta$:
$$\theta^*(H^p)=\overline{Span}_{H^p} \left[
\left\{\frac{\theta(\lambda)-\theta(z)}{\lambda-z}
\right\}_{\lambda\in\Bbb D}\right].$$
In the case $p=2$ we have $\theta^*(H^2)=H^2\circleddash \theta H^2$.

The following operator $V_\alpha$ acting from $L^2(\mu_\alpha)$ to
$\theta^*(H^2)$ was studied in [A2], [C] and [P1]:
$$V_\alpha f=\frac{Kf\mu_\alpha}{K\mu_\alpha}.$$
We will need the following 

\proclaim{Theorem 5 [C]}
For each $\alpha\in\Bbb T$
operator $V_\alpha$ maps $L^2(\mu_\alpha)$ onto
$\theta^*(H^2)$ unitarily. Operator $V_\alpha$ sends
operator $Y_\alpha$ of multiplication by $z$ in
$L^2(\mu_\alpha)$ into the following operator $T_\alpha$ in $\theta^*(H^2)$:
$$T_\alpha=
V_\alpha Y_\alpha V^*_\alpha f=z\left(f-\left(f,\frac\theta z\right)\frac\theta z\right)
+\left(f,\frac\theta z\right)\alpha.\tag 17$$
\endproclaim

Operators $T_\alpha$ are all possible unitary rank one perturbations of
the model contraction $T_\theta=P_\theta S$, where $S$ is the operator of
multiplication by $z$ in $\theta^*(H^2)$ and $P_\theta$ is the orthogonal
projector from $H^2$ to $\theta^*(H^2)$ (see [C]). Also note, that in
our situation $\theta(0)=0$, otherwise (17) would not be valid.

The conjugate operator $V_\alpha^*$ in (17)
can also be defined quite naturally:

\proclaim{Theorem 6 [P1]}
For each function $f\in\theta^*(H^2)$ and each $\alpha\in\Bbb T$
the non-tangential boundary values
of $f$ exist $\mu_\alpha$-a.e.\ Function $V_\alpha^*f\in L^2(\mu_\alpha)$
coincides with the non-tangential boundary values of $f$ $\mu_\alpha$-a.e.\
\endproclaim

Via the original Fourier transform and formula (17), each unitary
rank one perturbation $U_\alpha$ given by (16) is identified with an
operator in $\theta^*(H^2)$. Vector $\phi_1$ corresponds to
$1\in\theta^*(H^2)$ and vector $\phi_2$ corresponds to some function $f
\in\theta^*(H^2)$. As $\mu_\alpha$ was chosen to be the spectral measure
of $\phi_1$ for $U_\alpha$, the measure $|f|^2\mu_\alpha$ is
the spectral measure of $\phi_2$ for $U_\alpha$
(note that by Theorem 6 the non-tangential boundary values of
$|f|^2$ exist $\mu_\alpha$-a.e., therefore the notation $|f|^2\mu_\alpha$
makes sense).

Lemma 2 and Theorem 3 imply
that operator $U_{(\alpha,\beta)}=U_\alpha+\beta(\cdot,U_\alpha^{-1}
\phi_2)\phi_2$
is pure point for a.e.\ $\beta\in\Bbb T$ iff the (non-tangential
boundary values of the) derivative of $K|f|^2\mu_\alpha$ exist
a.e.\ on $\Bbb T$. We will use this fact in the proof of Theorem 4'.

By one of the properties of $\theta^*(H^2)$, if $f\in\theta^*(H^2)$
is equal to 0 at 0, then the function
$\hat f$ such that $\hat f=\theta \overline f$ a.e.\ on $\Bbb T$ also
belongs to $\theta^*(H^2)$.

We will need the following

\proclaim{Lemma 7} Let  $f\in\theta^*(H^2),||f||=1, f_0=f-f(0)$. Then

1)  there exist $g, h\in \cap_{0<p<1}\theta^*(H^p)$
such that $f_0\hat f_0=g+\theta h$;

2) $\hat f_0$ is equal to $\alpha\overline f_0$ $\mu_\alpha$-a.e.\
for every $\alpha\in\Bbb T$;

3) 
$$K|f|^2\mu_\alpha=\frac{g+\alpha h+f(0)\hat f_0+
\alpha\overline {f(0)}f_0+\alpha|f(0)|^2}
{\alpha-\theta}\tag 18$$
for every $\alpha\in\Bbb T$
\endproclaim
\remark{\bf Remark}
Recall, that $\theta(0)=0$, which implies $f_0\in \theta^*(H^2)$.
\endremark
\demo{Proof}
Since $f_0\hat f_0\in\theta^*(H^1)$, 1) follows from Riesz' Theorem.

To prove 2) it is enough to show that if $f_0$ is real $\mu_1$-a.e.,
then $\hat f_0=f_0$. But if $f_0$ is real $\mu_1$-a.e.\ then 
$$\theta\overline{\left(\frac{Kf_0\mu_1}{K\mu_1}\right)}=
\frac12\frac\theta{1-\overline\theta}(Pf_0\mu_1-i Im Kf_0\mu_1)=$$
$$=\frac12\frac1{1-\theta}(Pf_0\mu_1-i Im Kf_0\mu_1)=
-\frac{Pf_0\mu_1-i Im Kf_0\mu_1}{K\mu_1}$$
a.e.\ on $\Bbb T$.
Since $f_0\mu_1$ is a singular measure, Poisson integral $Pf_0\mu_1$
is equal to 0 a.e.\ on $\Bbb T$. Therefore
$$\hat f_0=\theta\overline{\left(\frac{Kf_0\mu_1}{K\mu_1}\right)}=
-\frac{Pf_0\mu_1-i Im Kf_0\mu_1}{K\mu_1}=
\frac{Kf_0\mu_1}{K\mu_1}=f_0$$
a.e.\ on $\Bbb T$.

To prove 3), let us first assume that $f\in H^\infty$. Then $g$, $h$ and
$g+\alpha h$ belong to $\theta^*(H^2)$. Since $\hat f_0=\alpha
\overline {f_0}$
 $\mu_\alpha$-a.e., $(f_0+f(0))(\hat f_0+\alpha\overline{f(0)})
 =\alpha |f|^2$ $\mu_\alpha$-a.e.\ Since
 $\theta=\alpha$ $\mu_\alpha$-a.e., $f_0\hat f_0=g+\theta h=g+\alpha h$
 $\mu_\alpha$-a.e.\ Hence
 $$g+\alpha h+f(0)\hat f_0+\alpha\overline {f(0)}f_0+
 \alpha|f(0)|^2=\alpha |f|^2$$
  $\mu_\alpha$-a.e.\
Thus, by Theorem 6,
$$\frac{K|f|^2\mu_\alpha}{K\mu_\alpha}=\overline\alpha(
g+\alpha h+f(0)\hat f_0+\alpha\overline {f(0)}f_0+\alpha|f(0)|^2).$$
Since $K\mu_\alpha=\frac 1{1-\overline\alpha\theta}$ (recall that $\theta(0)
=0$), we obtain (18).

In the general case, consider $f_n\to f$ in $\theta^*(H^2)$ where
$f_n\in H^\infty$. Then for each $f_n$ formula (18) will hold with
some functions $g_n$ and $h_n$ on the right hand side such that
$g_n+\alpha h_n\to g+\alpha h$ pointwise in $\Bbb D$. Since by Theorem 5
$|f_n|^2\to|f|^2$
in $L^1(\mu_\alpha)$,
$$\frac{K|f_n|^2\mu_\alpha}{K\mu_\alpha}\to
\frac{K|f|^2\mu_\alpha}{K\mu_\alpha}$$
pointwise in $\Bbb D$.
 $\blacktriangle$.
\enddemo

We will only prove Theorem 4' and other results from this section
for the case $n=2$.
This restriction will allow us to significantly shorten the proofs  without
losing any essential ideas.

\demo{Proof of Theorem 4' for n=2} WLOG $||\phi_i||=1$.
Let $f\in\theta^*(H^2)$ be the function corresponding to $\phi_2$
(we keep the notation introduced before Lemma 7).
Denote by $\mu_\alpha$ and $\nu_\alpha$ the spectral measures
 of $\phi_1$ and $\phi_2$ for
$$U_\alpha=U+(\alpha-1)(\cdot,U^{-1}\phi_1)\phi_1$$
respectively. Then $\nu_\alpha=|f|^2\mu_\alpha$.
To simplify (18) and other formulas we will assume that $f(0)=0$.
The general case can be treated in the same way.

If $f(0)=0$ then
Lemma 7 implies that
there exist $g, h\in \cap_{0<p<1}\theta^*(H^p)$ such that $f\hat f=
g+\theta h$ and
for every $\alpha\in\Bbb T$
$$K\nu_\alpha=\frac{g+\alpha h}{\alpha-\theta}.\tag 19$$
Let $\nu_{(\alpha,\beta)}$ be the spectral measure of $\phi_2$ for
$$U_{(\alpha,\beta)}=U_\alpha+(\beta-1)(\cdot,U_\alpha^{-1}\phi_2)\phi_2.$$
Using formulas (10) and (19) one can conclude that
$$K\nu_{(\alpha,\beta)}=\frac{\beta\frac{g+\alpha h}{\alpha-\theta}}
{1+(\beta-1)\frac{g+\alpha h}{\alpha-\theta}}.\tag 20$$

$1)\Rightarrow 2).$
As we discussed earlier, condition 1) means that the derivative
of the function on the right hand side of (19) exists a.e.\ on $\Bbb T$
for a.e.\ $\alpha\in\Bbb T$. Let us first show that
this implies that the derivatives of
$f$, $g$ and $\theta$ exist a.e.\ on $\Bbb T$.

Indeed,
$$\left(\frac{g+\alpha h}{\alpha-\theta}\right)'=
\frac 1{\alpha-\theta}g'+\frac \alpha{\alpha-\theta}h'+
+\frac{g+\alpha h}{(\alpha-\theta)^2}\theta'.\tag 21$$
Choose $\xi\in\Bbb T$ such that there exist finite non-tangential limits
$g(\xi)=a$, $h(\xi)=b$, $\theta(\xi)=c$ and
$\left(\frac{g+\alpha h}{\alpha-\theta}\right)'(\xi)$ for a.e.\ $\alpha$
(note that for a.e.\ $\xi$ such limits
exist for a.e.\ $\alpha$). Suppose that at least one of the non-tangential
limits $g'(\xi),\ h'(\xi)$ or $\theta'(\xi)$ does not exist.
WLOG we can assume that there exists a
sequence $\{z_n\}$ tending to $\xi$ non-tangentially,
such that the limit
$$\lim_{n\to\infty}g'(z_n)\tag 22$$
does not exist and both sequences
$\{h'(z_n)\}$ and $\{\theta'(z_n)\}$ are $O(g'(z_n))$ (such sequence
$\{z_n\}$ exists for either $g'$, $h'$ or $\theta'$).

Consider the vector function
$$v(\alpha)=\left(\frac1{\alpha-c},\frac\alpha{\alpha-c},
\frac{a+\alpha b}{(\alpha-c)^2}\right)$$
whose coordinate functions represent the limits of the coefficients
on the right hand side of (21) at $\xi$. Pick $\alpha_1,\alpha_2$ and
$\alpha_3$ such that the vectors $v(\alpha_i)$ are linearly independent
and the limit
$$\lim_{n\to\infty}
\left(\frac{g+\alpha_i h}{\alpha_i -\theta}\right)'(z_n)\tag 23$$
exists for $i=1,2,3$. Let constants $d_1,d_2,d_3$ be such that
$\sum_{i=1}^3 d_iv(\alpha_i)=(1,0,0)$. 
Then the sum
$\sum_{i=1}^3 d_i
\left(\frac{g+\alpha_i h}{\alpha_i -\theta}\right)'$
can be represented as
$$\sum_{i=1}^3 d_i
\left(\frac{g+\alpha_i h}{\alpha_i -\theta}\right)'(z)=
k_1(z)g'(z)+k_2(z)h'(z)+k_3(z)\theta'(z)$$
where $k_1(z_n)\to 1$, $k_2(z_n)=o(1)$ and $k_3(z_n)=o(1)$.
Since the limit (21) does not exist and $h'(z_n), \theta(z_n)=O(g'(z_n))$,
the limit
$$\lim_{n\to\infty}\sum_{i=1}^3 d_i
\left(\frac{g+\alpha_i h}{\alpha_i -\theta}\right)'(z_n)$$
does not exist. But by the choice of $\alpha_i$ the limit (23) exists
for $i=1,2,3$ and we have a contradiction. Hence the non-tangential
limits $g'(\xi)$, $h'(\xi)$ and $\theta'(\xi)$ exist for a.e. $\xi\in\Bbb T$.

Since
$$K\mu_1=\frac 1 {1-\theta},\tag 24$$
the derivative $(K\mu_1)'$ exists a.e.\ on $\Bbb T$.
Therefore
$U_\alpha=U+\alpha(\cdot,\phi_1)\phi_1$
is pure point for a.e.\ $\alpha\in\Bbb T$.
Also, since
$$K\nu_1=\frac{g+h}{1-\theta},\tag 25$$
the derivative $(K\nu_1)'$ exists a.e.\ on $\Bbb T$. Hence
$U_\alpha=U+\alpha(\cdot,\phi_2)\phi_2$
is pure point for a.e.\ $\alpha\in\Bbb T$.

$2)\Rightarrow 1)$
Since $U_\alpha=U+\alpha(\cdot,\phi_1)\phi_1$
is pure point for a.e.\ $\alpha\in\Bbb T$, 
the derivative $(K\mu_1)'$ exists a.e.\ on $\Bbb T$.
Hence by (24), $\theta'$ exist a.e. on $\Bbb T$.
Also
since $U_\alpha=U+\alpha(\cdot,\phi_2)\phi_2$
is pure point for a.e.\ $\alpha\in\Bbb T$, 
the derivative $(K\mu_2)'$ exists a.e.\ on $\Bbb T$.
Together with (25) this implies that $(g+h)'$ exists 
a.e.\ on $\Bbb T$.

The condition that
the derivatives $(K\mu_1)'$ and
$(K\nu_1)'$ exist a.e.\ on $\Bbb T$ is equivalent to
$$\int_\Bbb T \frac {d(\mu_1+\nu_1)(\xi)}{|\psi-\xi|^2}<\infty\tag 26$$
for a.e.\ $\psi\in\Bbb T$. Since $\mu_1+\nu_1=(1+|f|^2)\mu_1$
and $|f|<1+|f|^2$, (26) implies
$$\int_\Bbb T \frac {|f|d\mu_1(\xi)}{|\psi-\xi|^2}<\infty.$$
Therefore the derivatives $(Kf\mu_1)'$  and $(K\overline f\mu_1)'$
exist a.e.\ on $\Bbb T$. Since $f=(1-\theta)Kf\mu_1$ and
$\hat f=(1-\theta)K\overline f\mu_1$, $f'$ and $\hat f'$ also exist a.e.\ on $\Bbb T$.
Thus $(\hat f f)'=(g+\theta h)'$ exist a.e.\ on $\Bbb T$.
Since $(g+h)'$ and $(g+\theta h)'$ exist a.e.\ on $\Bbb T$,
$((1-\theta)h)'$ exists a.e.\ on $\Bbb T$. Since $\theta'$ and
$(g+h)'$ exist a.e.\ on $\Bbb T$, this implies that $g'$ and $h'$
exist a.e.\ on $\Bbb T$.
Therefore by (19), for each $\alpha\in\Bbb T$
 $(K\nu_\alpha)'$ exists a.e.\ on $\Bbb T$.
 By Lemma 2
this means that the operator
$U_{(\alpha,\beta)}$ is pure point for a.e.\ $\beta\in\Bbb T$ for every
$\alpha\in\Bbb T$.
$\blacktriangle$\enddemo

Theorem 4 implies that if operator $A_\lambda$ is pure point for
a.e.\ $\lambda\in \Bbb R^n$, then for any line $L\subset \Bbb R^n$ parallel
to one of the coordinate axis, $A_\lambda$ is pure point 
for a.e.\ $\lambda\in L$. We will now show that this statement holds
true even if one replaces the line $L$ with an arbitrary analytic curve.
\definition{Definition}
Let $I_1, I_2,...,I_n$ be inner functions in
the unit disk $\Bbb D$.
Let $\Sigma\in\Bbb T,\ m(\Sigma)=1$ be a set such that $I_k(\xi)$
exists for every $\xi\in\Sigma$ for any $1\leq k\leq n$.
Let $\gamma$ be a subset of the n-dimensional torus $\Bbb T^n$
such that
$$\gamma=\{\alpha\in\Bbb T^n | \alpha=\gamma(\xi)=(I_1(\xi), I_2(\xi),
...,I_n(\xi)),\ \xi\in\Sigma\}.$$
We will call such $\gamma$ an analytic curve in $\Bbb T^n$.

Similarly, if $J_1, J_2,...,J_n$ are analytic functions in $\Bbb C_+$
whose boundary values are real a.~e. on $\Bbb R$ and imaginary
parts are non-negative in $\Bbb C_+$, then we can consider
$\gamma\in\Bbb R^n$ such that
$$\gamma=\{\lambda\in\Bbb R^n | \lambda=\gamma(x)=(J_1(x), J_2(x),
...,J_n(x)),\ x\in\Sigma\}$$
where $\Sigma\in\Bbb R$ is a set of full measure such that $J_k(x)$
exists for every $x\in\Sigma$ for any $1\geq k\geq n$. We will
call such $\gamma$ an analytic curve in $\Bbb R^n$.


\enddefinition

Let $\omega$ be the standard conformal map from $\Bbb C_+$ to $\Bbb D$:
$\omega(z)=\frac{z-i}{z+i}$.
As one can see, for any analytic curve $\gamma\in\Bbb R^n, \gamma(t)=
(J_1(t),...,J_n(t))$
there exists an analytic curve $\eta\in\Bbb T^n, \eta(\xi)=
(I_1(\xi),...,I_n(\xi))$ such
that $J_k(z)=\omega^{-1} I_k(\omega(z))$. Conversely, any analytic
curve in $\Bbb T^n$ is an ``image'' of an analytic curve in $\Bbb R^n$.

We will need the following  generalization of formula (9):

\proclaim {Lemma 8}
Let $U, \phi_i (i=1,...,n)$ and $U_\alpha$ be as in Theorem 4'.
Denote by $\nu_\alpha$ the spectral measure of $\phi_n$ for
$U_\alpha$. Let $\gamma=(I_1,I_2,...,I_n)$ be an analytic curve in
$\Bbb T^n$. Then there exists a bounded analytic in $\Bbb D$ function $\varphi$
such that
$$\int_{\Bbb T} \nu_{\gamma(\xi)} dm(\xi)=\varphi m$$
i. e. for any Borel set $B\subset\Bbb T$
$$\int_{\Bbb T} \nu_{\gamma(\xi)}(B) dm(\xi)=\int_B\varphi d m.\tag 27$$
If $I_k(0)=0$ for $k=1,...,n$ then $\varphi= 1$.
\endproclaim

In the proof
we will obtain an explicit formula for the function $\varphi$ for the case
$n=2$, see (32).

\demo{Proof for n=2}
WLOG $||\phi_i||=1$.
Starting as in the proof of Theorem 4' (see (20)) we can observe
that
$$K\nu_{(I_1(\xi),I_2(\xi))}=
\frac{\frac{\overline{I_1(\xi)}g+ h}{1-\overline{I_1(\xi)}\theta}}
{\overline{I_2(\xi)}+(1-\overline{I_2(\xi)})
\frac{\overline {I_1(\xi)} g+ h}
{1-\overline{I_1(\xi)}\theta}}\tag 28$$
where $g, h\in \cap_{0<p<1}\theta^*(H^p)$ (to simplify
formula (18) we will again assume that the function from
$\theta^*(H^2)$ corresponding to $\phi_2$ is 0 at the origin).
If $\Cal K_z$ is Cauchy kernel for $z\in\Bbb D$, then
$$\int_\Bbb T\int_\Bbb T \Cal K_z(\omega) d\nu_{(I_1(\xi),I_2(\xi))}(\omega)dm(\xi)=
\int_\Bbb T
\frac{\overline{I_1(\xi)}\frac{g(z)+ h(z)}{1-\overline{I_1(\xi)}\theta(z)}}
{\overline{I_2(\xi)}+(1-\overline{I_2(\xi)})
\frac{\overline {I_1(\xi)} g(z)+ h(z)}
{1-\overline{I_1(\xi)}\theta(z)}}dm(\xi)
.\tag 29$$
Note that if we fix $\xi$ such that $|I_1(\xi)|=1$ then 
$$Re
\frac{\overline {I_1(\xi)} g(z)+ h(z)}
{1-\overline{I_1(\xi)}\theta(z)}>\frac 12\tag 30$$
 in $\Bbb D$ since fraction
$$\frac{\overline {I_1(\xi)} g(z)+ h(z)}
{1-\overline{I_1(\xi)}\theta(z)}\tag 31$$
is the Cauchy integral of a probability measure (namely it is the Cauchy integral
of the spectral
measure of $\phi_2$ for $U+I_1(\xi)(\cdot,U^{-1}\phi_1)\phi_1$, see
formula (19)).

Now if we fix $z\in\Bbb D$ then fraction (31)
 is a bounded antianalytic function of $\xi$ in $\Bbb D$.
Hence (30) must hold true for
every $\xi\in\Bbb D$ because it holds a.e.\ on $\Bbb T$. 
This implies that for any fixed $z\in\Bbb D$ the denominator
$$\overline{I_2(\xi)}+(1-\overline{I_2(\xi)})
\frac{\overline {I_1(\xi)} g(z)+ h(z)}
{1-\overline{I_1(\xi)}\theta(z)}$$
on the right hand side of (29) is an antianalytic function in $\Bbb D$ whose
absolute value is bounded away from 0. Hence for every $z\in\Bbb D$
the whole fraction
on the right hand side of (29) is a bounded antianalytic function
of $\xi$ in $\Bbb D$.
Thus
$$\int_\Bbb T\int_\Bbb T \Cal K_z(\omega) d\nu_{(I_1(\xi),I_2(\xi))}(\omega)dm(\xi)=
\frac{\overline{I_1(0)}\frac{g(z)+ h(z)}{1-\overline{I_1(0)}\theta(z)}}
{\overline{I_2(0)}+(1-\overline{I_2(0)})
\frac{\overline {I_1(0)} g(z)+ h(z)}
{1-\overline{I_1(0)}\theta(z)}}=\tag 32$$
$$=\varphi(z)=\int_{\Bbb T}
\Cal K_z\varphi dm.$$
Inequality (30) implies that $|\varphi|\leq\frac 1{1-|I_2(0)|}$.
Proceeding as in Examples 1 and 2 from the previous section,
we can replace $\Cal K_z$ in (32) with the characteristic function of
any Borel set $B\subset \Bbb T$ to obtain (26).
$\blacktriangle$\enddemo

We are now ready to generalize $1)\Rightarrow 2)$ part of Theorems 4
and 4' in the following way:

\proclaim{Theorem 9}
Let $A$ be a cyclic self-adjoint operator, $\phi_1,\phi_2,...,\phi_n$
its cyclic vectors. Let $\gamma$ be an analytic curve in $\Bbb R^n$.
Suppose that  operator
$$A_\lambda=A+\sum_{k=1}^n \lambda_k(\cdot,\phi_k)\phi_k$$
is pure point for a.e.\
$\lambda=(\lambda_1,\lambda_2,...,\lambda_n)
\in\Bbb R^n$. Then $A_{\gamma(t)}$ is pure point for a.e.\ $t\in\Bbb R$.
\endproclaim

\proclaim{Theorem 9'}
Let $U$, $\phi_1,\phi_2,...,\phi_n$ and $U_\alpha$ be the same
as in Theorem 4'. Let $\gamma$ be an analytic curve in $\Bbb T^n$.
Suppose that  operator $U_\alpha$
is pure point for a.e.\
$\alpha=(\alpha_1,\alpha_2,...,\alpha_n)
\in\Bbb T^n$. Then $U_{\gamma(\xi)}$ is pure point for a.e.\ $\xi\in\Bbb T$.
\endproclaim

\demo{Proof of Theorem 9' for n=2}
If $\nu_{(I_1(\xi),I_2(\xi))}$ is the spectral measure
of $\phi_2$ for $U_{(I_1(\xi),I_2(\xi))}$ then, as we showed in the
previous proof, its Cauchy integral must satisfy
$$K\nu_{(I_1(\xi),I_2(\xi))}=
\frac{\frac{\overline{I_1(\xi)}g+ h}{1-\overline{I_1(\xi)}\theta}}
{\overline{I_2(\xi)}+(1-\overline{I_2(\xi)})
\frac{\overline {I_1(\xi)} g+ h}
{1-\overline{I_1(\xi)}\theta}}\tag 33$$
By Lemma 2 to show that
$\nu_{(I_1(\xi),I_2(\xi))}$ is pure point it is enough to show that
function
$$\frac{\overline {I_1(\xi)} g+ h}
{1-\overline{I_1(\xi)}\theta}\tag 34$$
has a non-tangential derivative $\nu_{(I_1(\xi),I_2(\xi))}$-a.e.\

Denote by $E$  the subset of $\Bbb T$ where the non-tangential derivative
of (34) does not exist.
As was shown in the proof of Theorem 4', the condition that $U_\alpha$
is pure point for a.e.\ $\alpha$ implies that functions $f, g $ and
$\theta$ have nontangential derivatives a.e.\ on $\Bbb T$. Hence
$m(E)=0$. By Lemma 8 (formula (27)) this means that
$\nu_{(I_1(\xi),I_2(\xi))}(E)=0$ for a.e.\ $\xi\in\Bbb T$.
Therefore $\nu_{(I_1(\xi),I_2(\xi))}$ is pure point for
a.e.\ $\xi\in\Bbb T$.
$\blacktriangle$\enddemo

\enddocument